\documentclass[12pt]{amsart}
\usepackage{amsmath,amssymb,amscd,array}

\def\url#1{\expandafter\string\csname #1\endcsname}
\usepackage{url}

\usepackage{DLdef1}

\newcommand {\fby} {{\mathfrak{by}}}
\newcommand {\fer} {{\mathfrak{er}}}

\newcommand {\fdy} {{\mathfrak{dy}}}
\newcommand {\ffr} {{\mathfrak{fr}}}
\newcommand {\fmy} {{\mathfrak{my}}}

\newcommand {\fbr}{{\mathfrak{br}}}
\newcommand {\fbrj}{{\mathfrak{brj}}}

\newcommand{\un}{\underline{N}}
\newcommand{\del}{\partial}
\rmnameii{etr}{etr}
\rmnameii{evv}{ev}

\begin{document}
\title{Restricted simple Lie (super)algebras in characteristic $3$}

\author{Sofiane Bouarroudj${}^a$, Andrey Krutov${}^{b,e}$, Alexei Lebedev${}^c$, \\
Dimitry Leites${}^{a,d}$, Irina
Shchepochkina${}^e$}

\address{${}^a$New York University Abu Dhabi,
Division of Science and Mathematics, P.O. Box 129188, United Arab
Emirates; sofiane.bouarroudj@nyu.edu\\
$^b$Institute of Mathematics, Polish Aca\-demy of Sciences, ul. \'{S}niadeckich 8, 00-656 Warszawa, Poland;
a.o.krutov@gmail.com\\
${}^c$Equa
Simulation AB, R{\aa}sundav\"agen 100, Solna, Sweden; alexeylalexeyl@mail.ru\\
${}^{d}$De\-part\-ment of Mathematics, Stockholm University, SE-106 91
Stockholm, Sweden; dl146@nyu.edu\\
${}^e$Inde\-pendent University of Moscow, Bolshoj Vlasievsky per, dom
11, RU-119 002 Moscow, Russia; iparam53@gmail.com}



\begin{abstract} We  give explicit formulas proving restrictedness of  the following Lie (super)al\-ge\-bras: known
exceptional simple vectorial Lie (super)algebras in characteristic 3, 
deformed Lie (super)algebras with indecomposable Cartan matrix, and (under certain conditions) their simple subquotients over an algebraically closed field of characteristic 3, as well as one type of the deformed divergence-free Lie superalgebras with any number of indeterminates  in any characteristic.
\end{abstract}

\thanks{S.B. was partly supported by the grant AD 065 NYUAD. For the possibility to run difficult computations of this research  on the High Performance Computing Resources at NYUAD we are grateful to its Director, M.~Al Barwani. We thank S.~Skryabin and P.~Grozman for help}


\maketitle

\markboth{\itshape S. Bouarroudj\textup{,} A. Krutov\textup{,} A.
Lebedev\textup{,} D. Leites\textup{,} I.
Shchepochkina}{{\itshape Restricted Lie (super)algebras in characteristic $3$}}

\thispagestyle{empty}


\setcounter{section}{1}

\ssec{Introduction} Recall that a Lie aglebra $\fg$ over a field $\Kee$ of characteristic $p>0$ is called \textit{restricted} or \textit{with a $p$-structure} $x\mapsto x^{[p]}$, if for any $x\in\fg$, we have $(\ad_x)^p=\ad_{x^{[p]}}$ for some $x^{[p]}\in\fg$. A module $M$ over a restricted Lie algebra $\fg$, and representation $\rho$ defining $M$ are called \textit{restricted} if $\rho(x^{[p]})=\rho(x)^p$. A Lie superalgebra $\fg=\fg_\ev\oplus\fg_\od$ is \textit{restricted} if $\fg_\ev$ is restricted and $\fg$ is a restricted $\fg_\ev$-module. Thanks to \textbf{squaring}, i.e., the map $x\mapsto x^2(=\frac12[x,x]$  if $p\neq 2$) for any $x\in\fg_\od$, any restricted Lie superalgebra has a $2p$-structure, i.e., a map $x\mapsto x^{[2p]}$ for any $x\in\fg_\od$. 

In his Appendix to \cite{LL}, P.~Deligne advised us  to investigate first of all the \textbf{restricted} Lie (super)algebras and their \textbf{restricted} modules as related to geometry and hence of interest. This note is an addendum to  \cite{BLLS}, in which several general statements on restrictedness
valid for any $p>0$ are formulated, to \cite{BGL2},  where the Cartan matrices and Chevalley generators for modular Lie superalgebras are defined, and to \cite{BGL1, GL3, BGLLS, BGLLS1} describing Lie (super)algebras considered here. 
The  main result of \cite{BLLS} deals with $p=2$; here  we give  
examples for $p\neq 2$, mainly for $p=3$. 
The ground field is algebraically closed.

Classification \cite{BW} is implicit: to explicitly define $p|2p$-structure on a simple Lie superalgebra $\fg$ it suffices to give expressions of $w^{[p]}$ (resp. $w^{[2p]}$) for all even (resp. odd) elements of any basis of $\fg$. We give, at last, the explicit answer in case $\fsvect_{(1+\bar u)}(m;\One|2s)$, see~\eqref{new}; the deforms of series $\fh$ will be considered elsewhere. 

No classification of simple Lie superalgebras
is yet available for any $p>0$, or of simple Lie algebras  for $p=3$ and $2$, except for Lie (super)algebras with indecomposable Cartan matrix, and their simple subquotients, see \cite{BGL2}, whose $p|2p$-structure, if exists, is given explicitly. These Lie (super)algebras are \lq\lq symmetric\rq\rq, i.e., have a symmetric root system. For the classification of  \textit{true deforms}, i.e., \textit{results} of deformations which are neither trivial nor semitrivial,  of symmetric Lie (super)algebras whose restrictedness we establish here, see \cite{BLW, BGL3}. 

We also consider vectorial Lie (super)algebras. 
Following Bourbaki we use Gothic font for Lie (super)algebras;  $\One:=(1,\dots,1)$ is the shearing vector with the smallest heights of divided powers.
Proofs of lemmas, Fact~\eqref{Lie3}, formulas \eqref{Lmin1ps3str}  and \eqref{new}  
 are obtained with the help of the \textit{SuperLie} code, see \cite{Gr}. 

\ssec{Deforms of Lie (super)algebras with indecomposable Cartan matrix} In Lemmas \ref{br3} and \ref{g(1,6)}, the cocycles $c_k$ and the elements of the Chevalley basis $x_i$ (resp. $y_i$)  
corresponding to the positive (resp. negative) roots are given for the  Cartan matrix given in \cite{BGL3},
let the $h_j:=[x_j, y_j]$ be the elements of the maximal torus. 

\sssec{Deforms of $\fo(5)$ for $p=3$}\label{Exam}
Recall that the contact bracket of two divided powers $f,g\in \cO(p,q,t; \un)$
is defined to be
\begin{equation*}\label{cb}
 \{f,g\}_{k.b.}=\triangle f\cdot\partial_t g - \partial_t f\cdot\triangle g
 +\partial_p f\cdot\partial_q g - \partial_q f\cdot\partial_p g\text{~~with $\triangle f=2f - p\partial_p f - q\partial_q f$.}
\end{equation*}
A basis of $\mathfrak{L}(\eps, 0,0)$ is expressed in
terms of generating functions of $\mathfrak{k}(3; \One)$ and root vectors of $\fo(5)$ as
follows, see \cite[Prop.~3.2]{BLW}
\begin{equation*}\label{tab}\small
\renewcommand{\arraystretch}{1.4}
\begin{tabular}{|c|l|} \hline
$\deg$&the element${}_{\text{its weight}}\sim$ its generating function(=Chevalley basis vector) \\
\hline \hline
$-2$&$E_{-2\alpha -\beta}=[E_{-\alpha}, E_{-\alpha-\beta}]\sim 1(=y_4);$\\
\hline
$-1$& $E_{-\alpha}\sim p(=y_2);\quad E_{-\alpha-\beta}=[E_{-\beta}, E_{-\alpha}]\sim q(=y_3);$  \\
\hline 
$0$& $H_{\alpha}\sim -\eps t+ pq(=h_2);\quad H_{\beta}\sim -pq(=h_1); \quad
E_{\beta}\sim p^2(=y_1); \quad
E_{-\beta}\sim -q^2(=x_1);$  \\
\hline 
$1$& $E_{\alpha}\sim -(1+\eps)pq^2+ \eps
qt(=x_2);\quad E_{\alpha+\beta}=[E_\beta, E_\alpha]\sim (1+\eps)p^2q+\eps p t(=x_3);$  \\
\hline 
$2$& $E_{2\alpha+\beta}=[E_\alpha,
E_{\alpha+\beta}]\sim \eps(1+\eps)p^2q^2+ \eps^2t^2(=x_4).$  \\
\hline
\end{tabular}
\end{equation*}
\normalsize
Nonzero values of the deformed bracket with parameters $\delta$
and $\rho$ are as follows:
\begin{equation*}\label{croc1}
\renewcommand{\arraystretch}{1.4}\arraycolsep=2pt
\begin{array}{lll}
{}[E_{-2\alpha-\beta}, E_\beta] = \delta  E_{\alpha},
 &[E_{-2\alpha-\beta}, E_{-\alpha}]= \delta E_{-\beta},\quad{}
 &[E_{-\alpha}, E_{\beta}] = - \frac{\delta }{\eps} E_{2 \alpha+\beta} ,\\
 {}[E_{-2\alpha-\beta}, E_{-\alpha-\beta}]= \rho E_{\beta},\quad{}
 &[E_{-\alpha-\beta}, E_{-\beta}]= -\frac{\rho }{\eps}
 E_{2\alpha+\beta},
 &[E_{-2\alpha-\beta}, E_{-\beta}]= - \rho E_{\alpha+\beta}.
\end{array}
\end{equation*}
As proved in \cite{Kos}, in the family $\mathfrak{L}(\eps, \delta,\rho)$,  only $\mathfrak{L}(2,0,2)$ and Brown algebras $\fbr(2;\eps):=\mathfrak{L}(\eps,0,0)$ for $\eps\neq 0$ 
represent classes of non-isomorphic Lie algebras up to isomorphisms  
$\fbr(2;\eps)\simeq \fbr(2;\eps')$ if and
only if $\eps\eps'=1$ for $\eps\neq \eps'$; observe that $\fbr(2;-1)\simeq\fo(5)\simeq
\fsp(4)$.

\parbegin{Lemma}\label{3strLemm} 
The $3$-structure on $\mathfrak{L}(\eps, \delta,\rho)$  is given by the formulas
\begin{equation}\label{Lmin1ps3str}
\renewcommand{\arraystretch}{1.4}
\begin{array}{l}
h_1^{[3]}=h_1, \quad h_2^{[3]}=\eps^2 h_2, \ \
y_2^{[3]}= \delta(1+2 \eps^2)h_1+
\frac{\delta}{\eps}(1+2 \eps^2)h_2,\quad y_3^{[3]}=
\frac{\rho}{\eps}h_2, \\
 y_4^{[3]}=\eps \delta\rho(2+\eps^2)y_1,\quad
y_1^{[3]}= 
x_1^{[3]}=x_2^{[3]}=x_3^{[3]}=x_4^{[3]}=0.
\end{array}
\end{equation}
\end{Lemma}

In \cite{Kos}, Rudakov's claim \lq\lq $\mathfrak{L}(\eps, \delta,\rho)$ is restricted\rq\rq\ is cited but the explicit formulas \eqref{Lmin1ps3str} were never published, as far as we know. 

\sssbegin{Lemma}\label{br3} Let $\fg_{c_k}$ be the deform with even parameter $\lambda$ corresponding to the cocycle $c_k$ of $\fg=\fbr(3)$ or $\fbrj(2;3)$. The  Lie (super)algebras $\fg_{c_k}$, and those symmetric to them ($x\longleftrightarrow y$), are restricted. For any $k$, the $p|2p$-maps vanish on all weight vectors, except the following ones:  $h_i^{[3]}=h_i$ for all $i$ and   
 $x_3^{[3]}=- \lambda h_3$ for $\fbr(3)_{c_{-3}}$, and also $x_5^{[3]}=\lambda (h_2+h_3)$ for
$\fbr(3)_{c_{-6}}$; and $x_6^{[3]}=-\lambda (h_1+h_2+h_3)$ for $\fbr(3)_{c_{-9}}$; and $x_{10}^{[3]}=\lambda (h_1-h_2)$ for $\fbr(3)_{c_{-18}}$. Besides, for $\fbrj(2;3)_{c_{-12}}$, we have $x_{6}^{[3]}=2\lambda h_1$ and for $\fbrj(2;3)_{c_{-6}}$, we have $x_{3}^{[3]}=\lambda (h_1- h_2)$.
\end{Lemma}

\sssbegin{Lemma}\label{g(1,6)}  Let $\fg_{c_k}$ be the deform with odd parameter $\tau$ corresponding to the cocycle $c_k$ of $\fg=\fg(1,6)$ or $\fg(4,3)$ or $\fg=\fg(2,3)$. The Lie superalgebras $\fg_{c_k}$, corresponding to the cocycles $c_k$, and those symmetric to them ($x\longleftrightarrow y$),  are restricted. For any $k$, the $p|2p$-maps vanish on all weight vectors, except the following ones: 
$h_i^{[3]}=h_i$ for all $i$ for $\fg=\fg(1,6)$ and $\fg(4,3)$, and also $ d^{[3]}  ={}d$ for $\fg=\fg(2,3)$ modulo the central element ~$c=h_1-h_2$.
\end{Lemma}

\ssec{Fact}\label{fact} Let $\fg_0$ be a
restricted Lie (super)algebra and $\fg_{-1}$  an irreducible restricted $\fg_0$-module that generates the Lie superalgebra $\fg_{-} =\mathop{\oplus}_{-d\leq i<0}\fg_{i}$. Let vectorial Lie (super)algebra $\fg(\sdim;\un)$, where $\sdim$ is the superdimension of $\fg_-$, be the
\textit{prolong}, i.e., the result of generalized 
Cartan prolongation, see \cite{Shch}, of the pair $(\fg_-, \fg_0)$.  
It is easy to see that \textit{the Lie superalgebra
$\fg(\sdim;\un)$  is not restricted if} $\un\neq\underline{\One}$, see \cite{BLLS}; the proof of this statement for Lie algebras was first published in \cite[Th.2]{KfiD}.

\textbf{Fact}. \textit{If $\Zee$-graded vectorial Lie
(super)algebra $\fg:=\fg(\sdim;\One)$ --- the generalized Cartan prolong of its non-positive components, see \cite{Shch}, --- is restricted, and the $i$-th derived (super)algebra $\fg^{(i)}$
of $\fg$ contains a maximal torus of $\fg$, then $\fg^{(i)}$ is
restricted} \begin{equation}\label{pvert2p}
\begin{minipage}[l]{14cm}
(a) \textit{$\fh^{[p]}\subset \fh$; moreover, if the structure constants lie in
$\Zee/p$, then $h_i^{[p]}=h_i$ for
the basis elements $h_i$ of $\fh$};\\
(b) \textit{$w^{[p]}=0$ (resp. $w^{[2p]}=0$) for the other even (resp. odd)
weight elements $w$ of the basis of $\fg$ with weights relative a maximal torus of $\fder\ \fg$}.
\end{minipage}
\end{equation}
For the simple derived (super)algebra of every vectorial Lie (super)algebra $\fg$ we know for $p=3$, the $3$-structure is given by 
expressions~\eqref{pvert2p}, where $\fh$ is a maximal torus of $\fg_0$; for $p>3$, see \cite[Th.7.2.2]{S}.

\sssec{New examples}  The left column in \eqref{Lie3} shows where the simple  Lie (super)algebras in the right column are described for any $\un$; \textit{for them eqs.~\eqref{pvert2p} hold}:
\begin{equation}\label{Lie3}\footnotesize
\renewcommand{\arraystretch}{1.4}
\begin{tabular}{|l|l|}
\hline

\cite{GL3}&$\fdy^{(1)}(10; \One)$,\ $\fby^{(1)}(7; \One)$,\ $\fs\fby^{(1)}(7; \One)$,\ $\fmy(6; \One)$,\
$\fs\fmy^{(1)}(6; \One)$,\ $\fer^{(1)}(3; \One)$,\ $\ffr(3; \One)$\\

\hline
\cite{BL}, \cite{BGL1}, &$\mathfrak{Bj}(1; \One|7)$,\
$\mathfrak{Me}(3; \One|3)$,\ $\mathfrak{Brj}^{(2)}(3; \One|5)$,\ $\mathfrak{Bj}(4; \One|5)$,\ $\mathfrak{Bj}(3; \One|3)$,  \ $\mathfrak{Bj}^{(1)}(3; \One|4)$, \\
\cite{BGLLS}& $\mathfrak{Brj}^{(1)}(4|3)$ (only  $\un=\One$ is possible, unlike the above 2 lines)
\\
\hline
\end{tabular}
\end{equation}

\sssec{New formulas} 
The deforms of simple derived superalgebras of Lie  \textbf{super}algebras of the form $\fg(\sdim;\un)$ are not classified, we consider only one example. 
Let $\fsvect(m;\un|2s)$ be the Lie superalgebra  preserving the volume element $\vvol$
in even indeterminates $u_1, \dots,  u_m$ and odd ones $u_{m+1},\dots , u_{m+2s}$, let $\bar u=u_1^{(p^{\un_1}-1)}\cdots u_m^{(p^{\un_m}-1)}u_{m+1}\cdots u_{m+2s}$ and $\fsvect_{(1+\bar u)}(m;\un|2s)$ be the deform preserving $(1+\bar u)\vvol$. For $\fsvect_{(1+\bar u)}(m;\One|2s)$ and $\fsvect_{(1+\bar u)}^{(1)}(m;\One|2s)$ in characteristic $p>0$, we have (only distinctions with \eqref{pvert2p} are given): 
\be\label{new}
((1 - \bar{u})\del_i)^{[p] }= - (\del_i^{p-1} \bar{u}) \del_i, \text{~~where $\del_i$ is even.}
\ee


\label{lastpage}

\end{document}